\documentclass{amsart}
\usepackage{amssymb}

\newtheorem{theorem}{Theorem}[section]
\newtheorem{lemma}[theorem]{Lemma}
\theoremstyle{definition}
\newtheorem{definition}[theorem]{Definition}
\newtheorem{example}[theorem]{Example}

\newtheorem{corollary}[theorem]{Corollary}
\newtheorem{proposition}[theorem]{Proposition}

\numberwithin{equation}{section}
\newcommand{\beqa}{\begin{eqnarray*}}
\newcommand{\eeqa}{\end{eqnarray*}}
\newcommand{\beqn}{\begin{eqnarray}}
\newcommand{\eeqn}{\end{eqnarray}}

\newcommand{\pf}{\noindent {\it Proof.} }

\newcounter{cnt1}
\newcounter{cnt2}
\newcounter{cnt3}
\newcommand{\blr}{\begin{list}{$($\roman{cnt1}$)$}
        {\usecounter{cnt1} \setlength{\topsep}{0pt}
                \setlength{\itemsep}{0pt}}}
\newcommand{\bla}{\begin{list}{$($\alph{cnt2}$)$}
        {\usecounter{cnt2} \setlength{\topsep}{0pt}
                \setlength{\itemsep}{0pt}}}
\newcommand{\bln}{\begin{list}{$($\arabic{cnt3}$)$}
        {\usecounter{cnt3} \setlength{\topsep}{0pt}
                \setlength{\itemsep}{0pt}}}
\newcommand{\el}{\end{list}}
\newtheorem{thm}{Theorem}

\newtheorem{Def}[thm]{Definition}

\newtheorem{rem}[thm]{Remark}
\newcommand{\Rem}{\begin{rem} }
\newcommand{\bdfn}{\begin{Def} }
\newcommand{\edfn}{\end{Def}}

\title{Linear Hahn Banach Type Extension Operators in Banach Algebras of Operators}
\author{ Sudeshna Basu and Ajit Iqbal Singh }
\address{Department of Mathematics \\George Washington University\\ Washington DC 20052,USA \newline  Emeritus Scientist,  The Indian National Academy of Sciences,India }
\email{ sbasu@gwu.edu, sudeshna66@gmail.com,ajitis@gmail.com. }
\date{}
\begin{document}
\maketitle
\begin{abstract}
The notion of linear Hahn-Banach extension operator was first studied in detail by Heinrich and Mankiewicz (1982). Previously, J. Lindenstrauss (1966) studied similar versions of this notion in the context of non separable reflexive Banach spaces. Subsequently, Sims and Yost (1989) proved the existence of linear 
Hahn-Banach extension operators  via interspersing subspaces in a purely Banach space theoretic set up. In this paper, we study similar questions in the context of Banach modules and module homomorphisms, in particular, Banach algebras of operators on Banach spaces. Based on Dales, Kania, Kochanek, Kozmider and Laustsen(2013), and also Kania and Laustsen (2017), we give complete answers for reflexive Banach spaces and the non-reflexive space constructed by Kania and Laustsen from the celebrated Argyros-Haydon's space with few operators.
\end{abstract}
\section {Introduction}
Let $E$ be a Banach space over real or complex field and $E^*$, its dual. Let $B(E)$ denote the Banach algebra of all bounded linear operators on $E$. Let $\mathcal {A}$ denote a Banach algebra. Let $X$ denote a left Banach $\mathcal {A}$-module and $X'$, its ${\mathcal A}$ - dual, i.e. the space of all bounded left module homomorphisms from $X$ to $\mathcal {A}$ which is a right Banach $\mathcal {A}$-module. Let $\mathcal{L}(X)$ denote the space of all bounded $\mathcal {A}$-linear operators from X to itself.
%If $X$ is a Hilbert module, then $Adj(X)$ will denote all those bounded operators $ T \in B(X)$ which have an adjoint $T^*$.
We will freely use notation, terminology and basic results related to these notions from standard sources 
like,  \cite{DS}, \cite{La}, \cite{Pa}, \cite{D} and \cite{DDLS} . However, at times we will present some of them in the form that we need.

The following definition is well known in the literature, for reference see,\cite{SY1}.

\begin{definition} 
Let $E$ be a Banach space and let $M$ be a closed subspace of $E$. For each bounded linear functional $f \in M^*,$ we define 
$H_{M}(f) =  \{f'  \in E^* ~: \| f'\|= \|f\| , f'_{|M}= f.\}.$ Then it follows from Hahn-Banach theorem that, $H_{M}(f)$ is nonempty. It also follows that $H_{M}(f)$ is a $w^*$-compact and convex set.  Let $T : M^*\rightarrow X^*$ be such that $T(f)\in H_{M}(f)$. Then $T$ is called a {\it Hahn-Banach Extension operator.} 
\end{definition} 
Clearly such a $T$ is always norm preserving. It is natural to ask when some $T$ can be chosen to be linear. If $T$ is linear, it is called a { \it Linear Hahn Banach Extension operator.} 

We need the following definition:
\begin{definition}
By the density character, {\it dens X}, of a Banach space X, we mean the least
cardinality of any dense subset of X.
\end{definition}
We note the following result from \cite{SY1}.
\begin{theorem}\cite{SY1}\label {Sims and Yost}
Let $N$ be a closed subspace of a Banach space $E$. Then there exist a closed subspace M containing N, i.e. $N \subseteq M$ with  $dens~M = dens~N$ and a linear Hahn Banach extension operator 
$T:M^*\rightarrow E^*.$
\end{theorem}
\begin{rem}\label{interspersing}
In view of Theorem \ref{Sims and Yost}, we henceforth call the subspace $M$ as an { \it interspersing} subspace. 
\end{rem}
It is not difficult to show that if X is a Hilbert space, then T is linear. Conversely, if every subspace of  X  admits a linear  Hahn Banach Extension operator, then X is a Hilbert space. For details, see \cite{SY1}. For a Banach space linearity is not available in general but due to  Theorem \ref{Sims and Yost}, we can expect "plenty" of  subspaces which admit linear Hahn Banach Extension operators. 
 The notion of linear Hahn-Banach extension operator was first studied in detail by Heinrich and Mankiewicz \cite{HM}. Previously, J. Lindenstrauss studied similar versions of this notion in the context of non separable reflexive Banach spaces \cite{L}. Subsequently, Sims and Yost proved the existence of linear Hahn Banach extension operators  via interspersing subspaces in a purely Banach space theoretic set up \cite{SY1}. Sims and Yost also examined the existence of large families of
subspaces admitting linear Hahn-Banach extension operators, in duals of non-separable Asplund spaces \cite{SY2}. Harmand, Werner and Werner presented  linear Hahn-Banach extension operators in the context of M -ideals and M-embedded spaces \cite{HWW}. Much later, Lima and Oja studied the existence of linear Hahn Banach extension operators in the context of spaces  of bounded linear operators between two Banach spaces where one of them has the Radon Nikodym property \cite{LO}. Subsequently, Oja and Poldvere studied linear Hahn Banach extension operators  in the context of Principle of local reflexivity \cite{OP}. Daws studied similar extensions and Principle of Local Reflexivity for Banach modules and Banach algebras \cite{Da}.
In a recent article, Abrahamsen, studied the linear Hahn-Banach extension operators in the context of Almost Isometric Ideals in Banach spaces 
\cite{A}. Very recently, Yost has studied extensively, the notion of  linear Hahn-Banach extension operators and its relationship to ball intersection properties in the context of geometry of Banach spaces \cite{Y1},\cite{Y2}.

We have the following analogous definition for Banach modules.

\begin{definition}\label{module hahn banach}
Let $X$ be a Banach  left ${ \mathcal A}$ -module and $X'$ its ${ \mathcal A}$-dual.
\bla
 \item Suppose $M \subseteq X$ is a Banach left $\mathcal {A}$-submodule of $X$. Suppose an  $f \in M'$  has an extension,  
${ \tilde f} \in X'.$  We call ${ \tilde f}$, a { \it Hahn Banach Type Extension} of $f$. If $\|{ \tilde f} \| =\|f\|$, it is called  { \it Norm preserving Hahn Banach Type Extension} of $f$.
\item By the Hahn Banach extension set of f, we mean the set, $ {\mathcal H}_{M}(f) = \{ { \tilde f} \in X' ~:  { \tilde f}_{|M}= f\}.$ Suppose for each $f \in M'$ , ${\mathcal H}_{M}(f)$ is nonempty. Now, consider  $T : M'\rightarrow X'$ such that  $T(f)\in {\mathcal H}_{M}(f)$. Then $T$ is called a { \it Hahn-Banach Type Extension operator.} If $\|Tf\| =\|f\|$ for all $ f \in M ',$ then $T$ is called a   { \it Norm preserving Hahn-Banach Type Extension operator.}
\item For any Banach left ${\mathcal A}$ - submodule $N \subseteq M$, $M$ is said to  {\it intersperse} $N$ whenever $T:M'\rightarrow X'$ exists as a right 
$ \mathcal{A}$-linear map and $M$ is called an {\it interspersing left submodule.}
\el 
\end{definition}

\Rem
We can analogously define right and both sided, (norm preserving) Hahn Banach Type Extension, (norm preserving) Hahn-Banach Type Extension operator and interspersing Banach submodules.
\end{rem}

\begin{example}
A Hahn-Banach type extension may or may not exist in the context of Banach modules.
We give two prototypes to illustrate each case.
\blr
\item  Let $A$ be the C*-algebra of complex continuous functions on $[0,1]\cup[2,3].$ Then $A$ is a Hilbert  module over itself.
Let $M =\{f \in A :~ f = 0  \quad \mbox{on} \quad [2,3]\}.$ Then $M$ is a Hilbert $A$-submodule of $A$. Let $g= \chi_{[0,1]} \in M.$ For any $\phi \in M'$ we define $ \psi:A \rightarrow A$
as $\psi(u)= \phi (gu).$ Clearly $\psi$ is a Hahn-Banach Type Extension of $\phi.$ So we can now define  $T : M'\rightarrow A'$ as a Hahn Banach Type Extension operator as in the definition above.
\item Let $A= C[0,1]$ denote the C*-algebra of all complex continuous functions on $[0,1].$ Then $A$ is a Hilbert module over itself. Moreover, $M = C_{0}[0,1],$ the subset of all continuous functions vanishing at $0,$ is a Hilbert $A$-submodule
of $A$. For any $f \in M'$ we define a bounded  $A$-linear map $ \phi: M\rightarrow A$ given by  $\phi(f)(t)= f(t)sin(1/t), t \in (0,1], \phi(f)(0)=0, f \in M.$ Clearly $\phi$ cannot be extended to $A.$
\el
\end{example}

The following result is well known.
\begin {theorem}\label{standard}
Suppose $\mathcal{A}$ is a Banach  algebra and $X$ a left-Banach $\mathcal{A}$-module,then the set of bounded left module homomorphisms from $X$ to $\mathcal{A}$ is a right  $\mathcal{A}$-module via
$(f.a)(x)= f(x).a, \quad f \in  \mathcal {L}(X,\mathcal{A}), \quad x \in X,\quad  a \in \mathcal {A}.$
Also for $ g \in \mathcal {L}(\mathcal{A}, X)$, $f_{0}g $ is left module homomorphism of $\mathcal{A}$ to itself. Moreover, $ 
 f \rightarrow f_{0}g$ is a right module homomorphism of $\mathcal {L}(X,\mathcal{A})$ to $\mathcal {L}(\mathcal{A}, \mathcal{A}).$

\end{theorem}

The purpose of this paper is to investigate an analogue of  Theorem~\ref{Sims and Yost} in the context of Banach modules. We confine our discussion to Banach Algebra $\mathcal {A }$ of operators on a Banach space which is a module over itself and closed left ideals of $ \mathcal {A }$ i.e. closed left submodules.  As remarked above,  Hahn Banach type extensions for a module homomorphism may not be always possible, so the question  now boils down to three parts. First, we investigate the existence of Hahn-Banach type extensions. If there are such  extensions, then we investigate the existence of $ \mathcal{ A }$-linear Hahn Banach type extension operator, and we look for an  
{ \it interspersing} Banach $ \mathcal {A }$-submodule analogous to  Remark \ref {interspersing}. We address these results  in Section 2. Since  we cannot hope to have a result like   Theorem~\ref{Sims and Yost} in full generality, we  look at special cases in Section 3. 

Similar  results are studied  in various forms and in different contexts. For instance, see  \cite{AM}, \cite{B}, \cite{DDLS}, \cite{F}, \cite{FP}, \cite{H}, \cite{JS}, \cite{Li}, \cite{P}, and \cite{SS} and relevant references therein.

\section {Banach Algebra of Operators}
Let $E$ be a Banach space and $\mathcal{A}= B(E)$ be the Banach algebra of bounded linear operators from $E$ to itself.
For $x \in E, x^* \in E^*,$ the Banach space dual of $E,$ let $x^*\otimes x$ be the operator given by $x^*\otimes x(y)= x^*(y)x,$ for $y \in E.$ 
Then any $ T \in B(E)$ with finite dimensional range is given by $\sum_{j=1} ^{n}x^*_{j}\otimes x_j$ for some $ n \in N , x^*_{j} \in E^*, x_{j}\in E.$
We denote all these operators by $\mathcal{F}(E).$ Then $\mathcal{F}(E)$ is an ideal in $B(E).$ Fix any $x \in E,$ such that $\|x\| =1.$ By Hahn-Banach Theorem, there exists $\phi_{x} \in E^*$ such that $\|\phi_{x} \|=1$ and $\phi_{x} (x)=1.$  We fix any $\phi_{x}$ for the following discussion. We begin with  a  few concepts and results from \cite{DKKKL} that are basic and useful for our discussion.
Let $$M_{x}= \{ T\in B(E):~ T(x)=0\}$$ and $$N_{x}= \{ \phi_{x}\otimes y :~ y \in E\}.$$ Then $M_x$ and $N_x$ are closed linear subspaces of $\mathcal{A}$. 
In fact $M_x$ and $N_x$ are closed left ideals of  $\mathcal{A}.$  Also since  $\mathcal{A}$ is a Banach module over itself, $M_x$ and $N_x$ are left submodules of 
$\mathcal{A}.$ Also $M_{x} \neq 0$ if and only if  $dim E \geq 2.$ Further, $M_{x}\cap N_{x} = \{0\}.$

The following result is from \cite{DKKKL}.

\begin{theorem} \cite{DKKKL}\label{dales}
Let $\mathcal{A}=B(E) .$ Then for any maximal left ideal $\mathcal{I}$ of $\mathcal{A},$ exactly one of the following is true 
\begin{enumerate}
\item $\mathcal{I} =  M_{x}$ for some $x \in E, \|x\| =1.$
\item $\mathcal{I}$ contains  $\mathcal{F}(E).$
\end{enumerate}
\end{theorem}

We have the following theorem.
\begin{theorem}\label {banach algebra}
Let $ \mathcal{A}= B(E).$ Then for any $ \phi \in \mathcal {L}(M_{x}, \mathcal{A})$ there exists   $\tilde \phi
 \in \mathcal {L}(\mathcal{A}),$ and for any $ \psi \in \mathcal {L}(N_{x}, \mathcal{A})$ there exists 
 $\tilde\psi \in \mathcal{L}(\mathcal{A}, \mathcal{A})$ such that  $\tilde\phi_{|M_{x}}= \phi$ and $\tilde \psi_{|N_{x}}= \psi.$ Also $ \phi \rightarrow \tilde \phi$ and $\psi \rightarrow \tilde \psi$ are right module homomorphisms. In other words, $M_x$ and $N_x$ are interspersing left Banach submodules.
\end{theorem} 
\pf
Let $ T \in \mathcal{A}.$ Let $T_{2}= \phi_{x}\otimes Tx$ and $T_{1}= T- T_{2}.$
Then $T_{1} \in M_{x}$ and $ T_{2} \in N_{x}.$
Let $P:\mathcal{A}\rightarrow \mathcal{A} $ and $ Q:\mathcal{A}\rightarrow \mathcal{A}$ be defined by 
$P(T)= T_{1}$ and $Q(T)= T_{2}.$ Then P and Q are linear projections with range $M_{x}$ and $N_{x}$ respectively and $ I_{\mathcal{A}}= P+Q.$
For any $ T \in \mathcal{A} ,$ we have 
$$\|Q(T)\| = \|\phi_{x}\otimes Tx\|= \| \|\phi_{x}\|\|Tx\| \leq \|T\|.$$
Hence $\|Q\| \leq 1.$
But since $Q$ is a projection, $\|Q \| \geq 1.$
Thus we have $\|Q \|= 1.$
Consequently, $ 1 \leq \|P\| \leq 2.$
Now, for $T,S \in B(E),$
$Q(ST)= \phi_{x}\otimes(STx) = S(Q(T))$ as well.
Thus $P$ and $Q$ are left module homomorphisms of $ \mathcal{A}$ to itself.
Now suppose $\phi \in \mathcal {L}(M_{x}, \mathcal{A})$ and  $\psi \in \mathcal {L}(N_{x}, \mathcal{A}).$
Taking $\tilde \phi = \phi_{0}P$ and $\tilde \psi = \psi_{0}Q,$ we have the desired extensions from Theorem \ref{standard}.
\qed
\begin{rem}
\begin{enumerate}
\item  Clearly, $N_x$ being  a minimal left ideal and hence a minimal left  submodule, cannot act as an interspersing submodule for any other left submodule. 
\item $\|\psi\|= \|\tilde \psi\|$ in Theorem \ref{banach algebra} i.e., $\psi \rightarrow \tilde \psi $ is a norm preserving Hahn Banach Type Extension operator.
\end{enumerate}
\end{rem}

\begin{theorem} \label {singly generated}
Let  $L$ be a closed left ideal of $ \mathcal{A}= B(E) $ which is generated by a single  idempotent element in  $ \mathcal{A}.$
Then $L$ is an interspersing closed Banach left $ \mathcal{A}$-submodule of $ \mathcal{A}.$
\end{theorem}
\pf
Let $ D $ be an idempotent element in $ \mathcal{A}$ which generates $L$.
Then $D$ acts as a right identity on $L$.
For any $\phi \in L',$ we have $\phi (S)= \phi(SD)= S \phi(D)$ for all $S \in L.$  
Clearly $\phi$ extends to $\tilde {\phi} \in \mathcal{A}'.$
Hence the map $ \phi \rightarrow \tilde {\phi} $ is a Hahn Banach Type Extension Operator and $L$ is an interspersing Banach left $\mathcal{A}$-submodule.
\qed

We now prove the following lemma.

\begin {lemma}\label{finite ideal}
Let $E$ be a Banach space. Let $L$ be a closed left ideal of $\mathcal{A}= B(E)$ containing $ \mathcal{F}(E).$  Let $\phi$ be a bounded left module homomorphism of $L$ to $\mathcal{A}.$ If $\phi_{|\mathcal{F}(E)}= 0,$ then $ \phi=0$. Further, $\phi$ is extendible to $\tilde \phi$ if and only if $\phi_{|\mathcal{F}(E)}$ is extendible and $\tilde \phi$  is unique.
\end{lemma}

\pf
Indeed, let if possible, $\phi \neq 0,$ then  there exists $ T \in L$  with $\phi(T)\neq 0.$
So  there exists $u \in E$ such that $ \phi(T)(u) = \gamma  \mbox {(say)} \neq 0.$ Then there exists $f \in E^*$ such that $f(\gamma)= 1.$ 
Now,  $ f \otimes \gamma \in \mathcal{F}(E)$ and therefore $(f \otimes \gamma)T \in \mathcal{F}(E).$ 
So, $\phi((f\otimes \gamma)T)=0.$\\
But 
\begin{eqnarray*}
\phi ((f\otimes \gamma)T) &= & (f\otimes \gamma)\phi(T),\\
\mbox{so,} \quad \phi((f\otimes \gamma)T)(u)&= & (f\otimes \gamma)\phi(T)(u)\\
                           & = & f( \phi(T)(u))\gamma \\
                           & = &f(\gamma) \gamma \\
                           & = & \gamma \neq 0, \mbox{ a contradiction.}
\end{eqnarray*}
\qed

\begin{rem}\label {palmer-reflexive} 
\begin{enumerate}
\item Lemma \ref{finite ideal} sets the ground for adapting well-known material such as in Palmer (\cite{Pa}(1.7.14)), for ideals relevant to our discussion in the context of closed left ideals $L$ of $\mathcal {A} = B(E)$ containing $\mathcal{F}(E).$ We recall that the adjoint map $*$ from $\mathcal {A}$ to $B(E^*)$ is an anti homomorphism and it is surjective  if and only if $E$ is reflexive ( \cite{DS}, VI.9.13 ).
\item Let $\phi \in \mathcal {L}(L,\mathcal{A})$. Consider any $ g \in E^*$ and $ y \in E$ satisfying $g(y)=1.$ 
Then for $ f \in E^{*}, x \in E$,  $f \otimes x = (g \otimes x)(f \otimes y),$ so we have $$\phi(f \otimes x) = (g \otimes x)\phi (f \otimes y) =( g_{0}\phi)(f\otimes y)\otimes x = (\phi(f \otimes y))^{*}(g)\otimes x= h \otimes x \mbox {(say)}.$$ So the map $f \rightarrow h$ on $E^*$ to itself gives an $S \in B(E^*)$ that satisfies 
$\phi(f\otimes x) =S(f) \otimes x$ for $ x \in E$ and $f \in E^*.$ This, in turn, gives, $(\phi(T))^*= ST^*$ for $T \in \mathcal{F}(E).$    
\item Suppose $E$ is reflexive. Then $ S= V^*$ for some $ V \in \mathcal {A}.$ Consequently, $ (\phi(T))^* = V^*T^* = (TV)^*$ giving $\phi(T)= TV$ for $T \in \mathcal{F}(E).$ In view of Lemma \ref{finite ideal}, $\phi(T)= TV$ for $T \in L.$ Hence $ \phi \rightarrow V $ is an isometric antihomomorphism of $\mathcal {L}(L,\mathcal{A})$ to $\mathcal {A}.$ 
\item We now take $L= K(E)$, the closed ideal of compact operators in $\mathcal {A}.$ By \cite{DS} Chapter VI, for $W \in B(E^{**}), T \in L,  T^{**}(W)_{|E} \in L.$
For $S \in B(E^*)$, let $\phi_{S}(T)= T^{**}(S^*)_{|E}, T\in L.$ Then $\phi_{S} \in \mathcal {L}(L,\mathcal{A}).$ Further, $S \rightarrow \phi_{S}$ is an isometric isomorphism of $B(E^*)$ onto $\mathcal {L}(L,\mathcal{A}).$
\end{enumerate}
\end{rem}

\begin {proposition}\label{F(e)}
Let $E$ be a reflexive Banach space. Let $L$  be a closed left ideal of $\mathcal{A}= B(E)$ 
 such that   $\mathcal{F}(E) \subseteq L.$
Let $\phi$ be a bounded left module homomorphism of $L$ to $\mathcal{A}.$ Then there exists $\tilde \phi:\mathcal{A}\rightarrow \mathcal{A}$ such that 
$\tilde \phi$ is a norm preserving Hahn Banach type extension of $\phi.$ The map $\phi \rightarrow \tilde {\phi}$ from $\mathcal{L}( L, \mathcal {A})$ to 
$\mathcal{L}(\mathcal{A},\mathcal{A})$  is a right $\mathcal{A}$-linear norm preserving Hahn Banach type  Extension Operator. 
\end{proposition}

\pf 
We use Remark \ref{palmer-reflexive}(iii) above and set $S_{\phi}=V.$ Now define $\tilde { \phi}: \mathcal{A} \rightarrow \mathcal{A}$  by 
$\tilde {\phi}(T)=TS_{\phi}.$  Clearly $\tilde{\phi}$ is an extension of $\phi.$
The rest follows from Remark \ref{palmer-reflexive}.
\qed

In view of Theorem \ref{dales}, Theorem \ref{banach algebra} and Proposition \ref{F(e)}, we now have,

\begin {theorem} 
Let $E$ be a reflexive Banach space. Let $L$ be a closed left ideal of $\mathcal{A}= B(E), $ i.e. a  Banach left  $\mathcal{A}$ -submodule.
Then there exists an interspersing Banach left $\mathcal{A}$-submodule M  for $L$. 
\end{theorem}
\section { A special example}
In this section we look at a special Banach space, constructed in \cite{KL} from  Argyros and Haydon Banach space  $X_{AH}$ and a carefully chosen closed subspace $Y.$ 
$X_{AH}$ is a Banach space with few operators i.e. $B(X_{AH})= kI_{AH}+ K( X_{AH}),$ where the scalar fieled $k$ is the field of real or complex numbers and 
$K( X_{AH})$  is the space of compact operators from $X_{AH}$ to itself. For more details on $X_{AH}$ and $Y$ refer to \cite{KL}.
Let $Z = X_{AH}\oplus_{\infty}Y.$ Let $\mathcal {A} =B(Z).$ Then it follows from \cite{KL}, that every bounded operator $T \in {\mathcal A}$  has a unique representation 
\begin{displaymath} 
T=\left (\begin{array}{ccc}
\alpha_{11}I_{X_{AH}}+K_{11}&\alpha_{12}J+K_{12}\\
K_{21} & \alpha_{22}I_{Y}+K_{22}\end{array}\right) \ldots \ldots (*)
\end{displaymath}
where $\alpha_{11}, \alpha_{12}$ and $\alpha_{22}$ are scalars, $I_{X_{AH}}$ and $I_{Y}$ denote the identity operators on $X_{AH}$ and $Y$ respectively,
$J: Y \rightarrow X_{AH}$ is the inclusion map and the operators $K_{11}:X_{AH}\rightarrow X_{AH}$, $K_{12}:Y\rightarrow X_{AH}$, $K_{21}:X_{AH}\rightarrow Y$
and $K_{22}: Y\rightarrow Y$  are compact.

Let $M_{1}= \{T \in{\mathcal A}: \alpha_{22}=0\}$ and $M_{2}= \{ T \in {\mathcal A} :~ \alpha_{11}=0\}.$
Then $M_1$ and $M_2$ are closed  two sided ideals of codimension 1 in $B(Z).$

The following Theorem from \cite{KL} will be useful in our discussion .

\begin{theorem} \cite{KL}\label{kanya}The ideals $M_1$ and $M_2$ are the only maximal closed left ideals of $B(Z)= {\mathcal A}$ that contain $\mathcal{F}(Z).$
Also,
\begin {enumerate}
\item $M_1$ is generated as a left ideal by the two operators
\begin{displaymath}
A= \left ( \begin{array}{ccc}
I_{X_{AH}}&0\\
0 & 0\end{array}\right)
\end{displaymath}
and
\begin{displaymath}
B = \left (\begin{array}{ccc}
0&J\\
0 & 0\end{array}\right)
\end{displaymath}
but $M_1$ is not generated as a left ideal by a single bounded operator on $Z,$
\item  $ M_2$ is not finitely generated as a left ideal.
\end{enumerate}
\end{theorem}

\begin{lemma}\label{firstlemma}
Let $\phi:M_{1}\cap M_{2} \rightarrow { \mathcal A }$ be a left module homomorphism.
Then 
%\begin{displaymath}
%\phi(A)= \left ( \begin{array}{ccc}
%w_{11}& w_{12}\\
%0 & 0\end{array}\right)
%\end{displaymath}
%and

$\phi(B)$ has the form  $\left (\begin{array}{ccc}
v_{11}& v_{12}\\
0 & 0\end{array}\right).$
Also $Range( v_{11})\subseteq Y$ and $Range (v_{12})\subseteq Y.$
\end{lemma}

\pf Let 
%\begin{displaymath}
%\phi(A)= \left ( \begin{array}{ccc}
%w_{11}& w_{12}\\
%w_{21} & w_{22}\end{array}\right)
%\end{displaymath}
%and
\begin{displaymath}
\phi(B) = \left (\begin{array}{ccc}
v_{11}& v_{12}\\
v_{21} & v_{22}\end{array}\right).
\end{displaymath}

Since $B^2=0,$ 
\begin{displaymath}
B\phi(B) = \left (\begin{array}{ccc}
0& 0\\
0 & 0\end{array}\right).
\end{displaymath}
Also \begin{displaymath}
B\phi(B) =  = \left (\begin{array}{ccc}
0& J\\
0 & 0\end{array}\right)
\left (\begin{array}{ccc}
v_{11}& v_{12}\\
v_{21} & v_{22}\end{array}\right)
= \left (\begin{array}{ccc}
Jv_{21}& Jv_{22}\\
0 & 0\end{array}\right).
\end{displaymath}
So $Jv_{21}= 0$ which, in turn gives $v_{21}=0 ;$
and $Jv_{22}= 0$ which, in turn gives $v_{22}=0.$
Therefore \begin{displaymath}
\phi(B) =  
\left (\begin{array}{ccc}
v_{11}& v_{12}\\
0 & 0\end{array}\right).
\end{displaymath}
Now for any $T\in B(Z),$ 
\begin{displaymath}
TB = 
\left (\begin{array}{ccc}
T_{11}& T_{12}\\
T_{21} & T_{22}\end{array}\right)
\left (\begin{array}{ccc}
0& J\\
0 & 0\end{array}\right)
= \left (\begin{array}{ccc}
0& T_{11}J\\
0 & T_{21}J\end{array}\right).
\end{displaymath}

 Also, \begin{displaymath}
 T\phi(B) = \left (\begin{array}{ccc}
T_{11}& T_{12}\\
T_{21} & T_{22}\end{array}\right)
 \left (\begin{array}{ccc}
v_{11}& v_{12}\\
0 & 0\end{array}\right)
= \left (\begin{array}{ccc}
T_{11}v_{11}& T_{11}v_{12}\\
T_{21}v_{11} & T_{21}v_{12}\end{array}\right).
\end{displaymath}

Let, if possible, $Range(v_{11})$ be not contained in $Y.$
Then there exists $x \in X_{AH}$ such that $s= v_{11}(x)\notin Y.$
Hence there exists $f \in X^{*}_{AH}$ such that $f(s)\neq 0$ and $f(Y)=\{0\}.$
Fix any $z \in X_{AH}$, $z \neq 0.$
Let $T_{11}= f\otimes z.$
For any $ y \in Y,$
$T_{11}J(y) = (f\otimes z)(y) = f(y)z =0 .$ 

Let $T_{21}= 0,$
then 
\begin{displaymath}
TB = \left (\begin{array}{ccc}
0& 0\\
0 & 0\end{array}\right).
\end{displaymath}
Since $\phi(TB)= T\phi(B),$ we arrive at a contradiction.
So $Range(v_{11})\subseteq Y.$
Similarly  $Range(v_{12})\subseteq Y.$
\qed
\begin{lemma}\label{second lemma}
Let $\phi:M_{1} \rightarrow { \mathcal A }$ be a left module homomorphism.
Then 

$\phi(A)$ has the form $\left ( \begin{array}{ccc}
w_{11}& w_{12}\\
0 & 0\end{array}\right).
$
\end{lemma}

\pf Let 
\begin{displaymath}
\phi(A)= \left ( \begin{array}{ccc}
w_{11}& w_{12}\\
w_{21} & w_{22}\end{array}\right).
\end{displaymath}

Since $A^{2}=A,$ it follows that $A \phi(A) = \phi(A).$
Now 
\begin{displaymath}
  A\phi(A)= \left ( \begin{array}{ccc}
I_{X_{AH}}& 0\\
0 & 0\end{array}\right) \left ( \begin{array}{ccc}
w_{11}& w_{12}\\
w_{21} & w_{22}\end{array}\right)
=\left ( \begin{array}{ccc}
w_{11}& w_{12}\\
0& 0\end{array}\right).
\end{displaymath}
Hence $w_{21}=0$ and $w_{22}=0$ and
\begin{displaymath}
\phi(A)= \left ( \begin{array}{ccc}
w_{11}& w_{12}\\
0 & 0\end{array}\right).
\end{displaymath}
\qed
\begin{lemma}\label{third lemma}
Let $\phi:M_{2} \rightarrow { \mathcal A }$ be a left module homomorphism.
If \begin{displaymath}
D = \left ( \begin{array}{ccc}
0& 0\\
0 & I_{Y}\end{array}\right),
\end{displaymath}
then $\phi(D)$ has the form $\left ( \begin{array}{ccc}
0& 0\\
u_{21} & u_{22}\end{array}\right).$
\end{lemma}
\pf
The proof is  similar to that of  Lemma \ref{second lemma}. {\qed

%\begin{corollary}
%If $ \phi:M_{1} \rightarrow B(Z)$ is as above, then $Jv_{11}= v_{11},Jv_{12}= v_{12}.$
%\end{corollary}
\begin{lemma}\label{fourth lemma}
Let $\phi:M_{1}\rightarrow B(Z)$ be a left module homomorphism. 
Let \begin{displaymath}
C = \left (\begin{array}{ccc}
w_{11}& w_{12}\\
J^{-1}v_{11}& J^{-1}v_{12}\end{array}\right),
\end{displaymath} where $J^{-1}$ is defined on $Y$ considered as a subspace of $ X_{AH}$ to $Y$ in its own right. Then $\phi(T)= TC$ for all $T \in M_{1}.$
\end{lemma}
\pf
Let $T\in M_{1}.$
Then by Theorem \ref{kanya}, it follows that $T= RA + SB $ for some $R, S \in \mathcal {A}.$
So \begin{eqnarray} \phi(T) &= & \phi(RA+SB)\\
&=& R \phi(A)+ S \phi(B)\\
& =& R\left (\begin{array}{ccc}
w_{11}& w_{12}\\
0& 0\end{array}\right) + S\left (\begin{array}{ccc}
v_{11}& v_{12}\\
0& 0\end{array}\right)\\
&=& \left (\begin{array}{ccc}
R_{11}& R_{12}\\
R_{21}& R_{22}\end{array}\right)\left (\begin{array}{ccc}
w_{11}& w_{12}\\
0& 0\end{array}\right)+\left (\begin{array}{ccc}
S_{11}& S_{12}\\
S_{21}& S_{22}\end{array}\right)\left (\begin{array}{ccc}
v_{11}& v_{12}\\
0& 0\end{array}\right)\\
&=& \left (\begin{array}{ccc}
R_{11}w_{11}& R_{11}w_{12}\\
R_{21}w_{11}& R_{21}w_{12}\end{array}\right)+\left (\begin{array}{ccc}
S_{11}v_{11}& S_{11}v_{12}\\
S_{21}v_{11}& S_{21}v_{12}\end{array}\right).
\end{eqnarray}
Now \begin{displaymath}
TC =(RA+SB)C = RAC + SBC.
\end{displaymath}
Also 
\begin{eqnarray} RAC &=
 =& \left (\begin{array}{ccc}
R_{11}& R_{12}\\
R_{21}& R_{22}\end{array}\right)\left (\begin{array}{ccc}
I_{X_{AH}}& 0\\
0& 0\end{array}\right)  \left (\begin{array}{ccc}
w_{11}& w_{12}\\
J^{-1}v_{11}& J^{-1}v_{12}\end{array}\right)\\
&=& \left (\begin{array}{ccc}
R_{11}& 0\\
R_{21}& 0\end{array}\right)\left (\begin{array}{ccc}
w_{11}& w_{12}\\
J^{-1}v_{11}& J^{-1}v_{12}\end{array}\right)\\
&=& \left (\begin{array}{ccc}
R_{11}w_{11}& R_{11}w_{12}\\
R_{21}w_{11}& R_{21}w_{12}\end{array}\right).
\end{eqnarray}

Similarly, 
\begin{eqnarray} SBC & =
 =& \left (\begin{array}{ccc}
S_{11}& S_{12}\\
S_{21}& S_{22}\end{array}\right)\left (\begin{array}{ccc}
0 & J\\
0& 0\end{array}\right)  \left (\begin{array}{ccc}
w_{11}& w_{12}\\
J^{-1}v_{11}& J^{-1}v_{12}\end{array}\right)\\
&=& \left (\begin{array}{ccc}
S_{11}& S_{12}\\
S_{21}& S_{22}\end{array}\right)\left (\begin{array}{ccc}
v_{11}& v_{12}\\
0& 0\end{array}\right)\\
&=& \left (\begin{array}{ccc}
S_{11}v_{11}& S_{11}v_{12}\\
S_{21}v_{11}& S_{21}v_{12}\end{array}\right).
\end{eqnarray}

Adding (3.8) and (3.11) and equating to (3.5), it follows that 
$\phi(T)= TC.$
\qed

\begin{theorem}\label{M1}
Let $\phi:M_{1} \rightarrow B(Z)$ be any bounded left module homomorphism. Then there exists a right $\psi:B(Z)\rightarrow B(Z)$ which is a norm preserving Hahn Banach Type extension of $\phi.$ Also, there exists  $B(Z)$-linear norm preserving Hahn Banach Type Extension operator $ \xi :\mathcal{L}(M_{1},B(Z)) \rightarrow  \mathcal{L}(B(Z))$
given by $\xi(\phi)= \psi.$ Finally,  $M_1$ is an interspersing left Banach $ \mathcal A$-submodule.
\end{theorem}
\pf
From Lemma \ref{fourth lemma}, it follows that we can take, $ \psi(T)= TC$, $T  \in {\mathcal A}.$ The rest follows from Remark \ref{palmer-reflexive} and Lemma \ref{finite ideal}.\qed

For the following discussion, we introduce some notation.

Let $ \mathcal {A}_{1}= B(X_{AH})$, $ \mathcal {A}_{2}= B(Y,X_{AH})$, $ \mathcal {A}_{3}= K(X_{AH})$,  $ \mathcal {A}_{4}= B(X_{AH},Y)$, $ \mathcal {A}_{5}= K(X_{AH},Y)$ and $ \mathcal {A}_{6} =B(Y).$ Then $ \mathcal {A}_{2}$ and $\mathcal {A}_{3}$ are Banach left modules over the Banach algebra 
$\mathcal {A}_{1}.$ Also $\mathcal{A}_{4}$ and $\mathcal{A}_{5}$ are Banach left modules over $ \mathcal{A}_{6}.$ 
\begin{theorem} \label{m2}Let $\psi_{1}:\mathcal {A}_{3}\rightarrow \mathcal {A}_{1}$ and $\psi_{2}:\mathcal {A}_{3}\rightarrow\mathcal {A}_{2}$ be bounded left $\mathcal {A}_{1}$-module homomorphisms. Then, for $K_{21}\in\mathcal {A}_{5}$, Range $\psi_{1}(JK_{21})\subseteq Y$ and Range $\psi_{2}(JK_{21})\subseteq Y.$
For \begin {displaymath}
T= \left (\begin{array}{ccc}
T_{11}& T_{12}\\
T_{21}& T_{22}\end{array}\right) \in M_2\end{displaymath} define 
\begin {displaymath}
\phi(T)= \left (\begin{array}{ccc}
\psi_{1}(T_{11})& \psi_{2}(T_{11})\\
J^{-1}\psi_{1}(JT_{21})& J^{-1}\psi_{2}(JT_{21})\end{array}\right).\end{displaymath}
Then we have the following:
\blr
\item $\phi$ is a bounded left $\mathcal{A}$-module homomorphism of $M_2$ to $\mathcal{A}$.
\item 
$\phi$ extends to a bounded left $\mathcal{A}$-module homomorphism of $\mathcal{A}$ to itself if and only if $\psi_1$ and $\psi_2$ extend to bounded left $\mathcal{A}_{1}$-module homomorphisms of $\mathcal{A}_{1}$ to $\mathcal{A}_{1}$ and $\mathcal{A}_{2}$ respectively.
\el

On the other hand, if  $\phi: M_{2} \rightarrow \mathcal{A}$ is a bounded left module homomorphism, then there exist $U \in \mathcal{A}$ of the form 
\begin {displaymath}
U= \left (\begin{array}{ccc}
0& 0\\
u_{21}& u_{22}\end{array}\right)\end{displaymath} and bounded left $\mathcal {A}_{1}$-module homomorphisms, $\psi_{1}:\mathcal {A}_{3}\rightarrow \mathcal {A}_{1}$, $\psi_{2}:\mathcal {A}_{3}\rightarrow\mathcal {A}_{2},$ such that 
\begin {displaymath}
\phi(T)= \left (\begin{array}{ccc}
\psi_{1}(T_{11}) + T_{12}u_{21}& \psi_{2}(T_{11}) + T_{12}u_{22}\\
J^{-1}\psi_{1}(JT_{21}) + T_{22} u_{21}& J^{-1}\psi_{2}(JT_{21})+ T_{22}u_{22}\end{array}\right).\end{displaymath}.
%Where $\psi_{1}:\mathcal {A}_{3}\rightarrow \mathcal {A}_{1}$ and $\psi_{2}:\mathcal {A}_{3}\rightarrow\mathcal {A}_{2}$ be bounded left $\mathcal {A}_{1}$-module homomorphisms.

\end{theorem}
\pf
By \cite{KL}. Prop 4.1(iii)), $ K(Z)$ has  a bounded two-sided approximate identity consisting of finite rank canonical projections. 
We can apply  Remark 4(ii) to the bounded  map  $\psi_1$ and obtain that for  $K \in \mathcal{A}_3$ with $range K \subseteq Y,$ $ range(\psi_{1}(K))\subseteq Y.$ Similar arguments  will give the result for $\psi_{2}(K)$. We can now take $K = J K_{21}$ for $K_{21} \in \mathcal{A}_{5}.$

(i)Clearly $\phi$ is linear and bounded.
Also, for for any $ S \in \mathcal{A},$ 
 
\begin {displaymath}
SD= \left (\begin{array}{ccc}
0& S_{12}\\
0& S_{22}\end{array}\right)\end{displaymath} 
and therefore, $\phi(SD) =0.$ 
For $ S \in \mathcal{A}$, $T \in M_{2},$ we have,
\begin {displaymath}
ST= \left (\begin{array}{ccc}
S_{11}T_{11}+ S_{12}T_{21}& S_{11}T_{12}+ S_{12}T_{22}\\
S_{21}T_{11}+ S_{22}T_{21}& S_{21}T_{12}+ S_{22}T_{22}\end{array}\right).\end{displaymath} .

Similarly,
\begin {displaymath}
S\phi(T)= \left (\begin{array}{ccc}
S_{11}\psi_{1}(T_{11})+ S_{12}J^{-1}\psi_{1}( JT_{21}) & S_{11}\psi_{2}(T_{11})+ S_{12}J^{-1}\psi_{2}(JT_{21})\\
S_{21}\psi_{1}(T_{11}) + S_{22}J^{-1}\psi_{1}(JT_{21})& S_{21}\psi_{2}(T_{11})+ S_{22}J^{-1}\psi_{2}(JT_{21})\end{array}\right).\end{displaymath} .
As per the part (i) of the  proof of  Theorem 1.9 in [KL], for K compact on $Y$ to $X_{AH}$ or to $Y,$ there is a compact extension $\tilde K$ defined on  $X_{AH}$ to  $X_{AH}$ or to $Y$ respectively. Obviously $K = \tilde{K}J.$
      We  may express $S_{12}= a J+L_{12}$ with $L_{12}$ compact. Then $S_{12}T_{21}= a JT_{21}+ L_{12}T_{21}= (a I_{X_{AH}}+ \tilde{L_{12}})J T_{21}.$ Because $\psi_{1}$ is a left $\mathcal {A}_{1}$-module homomorphism, we have $\psi_{1}(S_{12} T_{21})= (aI_{X_{AH}}+ \tilde{L_{12}}) \psi_{1}(J T_{21})= (a I_{X_{AH}}+\tilde{L_{12}}) J J^{-1}(\psi_{1}(J T_{21}))=( a J+L_{12})  J^{-1} \psi_{1}(J T_{21}) =S_{12} J^{-1} \psi_{1}(J T_{21}). $
      To deal with $S_{22} J^{-1} \psi_{1}(J T_{21}),$ we first write $S_{22}= b I_{Y}+L_{22}$ with $ L_{22}$ compact. Then  
$J S_{22} T_{21}= J(b I_{Y} T_{21}+ L_{22} T_{21})=(b I_{X_{AH}}+ J \tilde{L}_{22}) J T_{21}.$ Because $\psi_{1} $is a left $\mathcal{A}_1$- module homomorphism, we have $\psi_{1}(J S_{22} T_{21}))= (b I_{X_{AH}}+J \tilde{L_{22}}) \psi_{1}(J T_{21})=(b I_{X_{AH}}+J \tilde{L_{22}}) J J^{-1} (\psi_{1}(J T_{21}))=( b J+J L_{22})  J^{-1}\psi_{1}(J T_{21}) =J S_{22} J^{-1} \psi_{1}(J T_{21}). $
   We can treat with other terms in a similar  manner.
So, 
\begin {eqnarray*}
S\phi(T)&=& \left (\begin{array}{ccc}
\psi_{1}(S_{11}T_{11})+ \psi_{1}(S_{12}T_{21}) & \psi_{2}(S_{11}T_{11})+ \psi_{2}(S_{12}T_{21})
\\
J^{-1}\psi_{1}(JS_{21}T_{11})+ J^{-1}\psi_{1}(JS_{22}T_{21})& J^{-1}\psi_{2}(JS_{21}T_{11})+ J^{-1}\psi_{2}(JS_{22}T_{21})\end{array}\right)\\
&=& \left (\begin{array}{ccc}
\psi_{1}(S_{11}T_{11} + S_{12}T_{21})& \psi_{2}(S_{11}T_{11}+ S_{12}T_{21})\\
J^{-1}\psi_{1}(J(S_{21}T_{11}+ S_{22}T_{21}))& J^{-1}\psi_{2}(J(S_{21}T_{11}+ S_{22}T_{21}))\end{array}\right)\\
&=&
\phi \left (\begin{array}{ccc}
S_{11}T_{11} +S_{12}T_{21}  & 0 \\
S_{21}T_{11}+ S_{22}T_{21}& 0\end{array}\right)\\
&=&\phi(ST).
\end{eqnarray*} .
Hence $\phi$ is a left $\mathcal {A}$-module homomorphism.

(ii) Suppose $\phi$ extends to a $\psi \in \mathcal{L(A)}.$ Then there exists $V \in \mathcal{A}$ such that $\psi(T)= TV$ for $T\in \mathcal{A}.$
So, for $T_{11} \in \mathcal{A}_{3}$
\begin {eqnarray*}
\left (\begin{array}{ccc}
\psi_{1}(T_{11}) & \psi_{2}(T_{11})\\
0& 0\end{array}\right)
&=& \phi (\left (\begin{array}{ccc}
T_{11}  & 0\\
0& 0\end{array}\right))\\
&=&
\left (\begin{array}{ccc}
T_{11}& 0\\
0& 0\end{array}\right)\left (\begin{array}{ccc}
V_{11}  & V_{12}\\
V_{21}  & V_{22}\\
\end{array}\right)\\
&=&\left (\begin{array}{ccc}
T_{11}V_{11} & T_{11} V_{12}\\
0 & 0\\
\end{array}\right).\\
\end{eqnarray*}
This gives us that $\psi_{1}(T_{11})= T_{11}V_{11}$ and $\psi_{2}(T_{11})= T_{11}V_{12}$ for $ T_{11} \in \mathcal{A}_{3}.$
Consequently $\psi_{1}$ and $\psi_{2}$ extend to $ \tilde {\psi_{1}}, \tilde {\psi_{2}}$ on $ \mathcal {A}_{1}$, given by 
$\tilde {\psi_{1}}(T_{11})= T_{11}V_{11}$, $\tilde{\psi_{2}}(T_{11})= T_{11}V_{12},$ for $ T_{11}\in \mathcal {A}_{1}.$

Conversely, suppose $\psi_{1}$ and $\psi_{2}$ extend to $ \tilde {\psi_{1}}\in \mathcal {L(A}_{1})$ and $ \tilde {\psi_{2}} \in \mathcal{L} (\mathcal {A}_{1},\mathcal {A}_{2}) $ respectively.

Set  $V_{11}= \tilde {\psi_{1}}(I_{X_{AH}})$, $V_{12} = \tilde {\psi_{2}}(I_{X_{AH}})$ and \begin {displaymath}
V= \left (\begin{array}{ccc}
V_{11}& V_{12}\\
0& 0\end{array}\right).\end{displaymath}.

Then $\tilde{\psi_{1}}(T_{11})= T_{11}V_{11}$ and $ \tilde {\psi_{2}}(T_{11})= T_{11}V_{12}$ for $ T_{11} \in \mathcal {A}_{1}.$
So, for $ T_{11}\in \mathcal {A}_{3}, T \in M_{2},$
\begin {eqnarray*}
\phi(T) &=&\left (\begin{array}{ccc}
T_{11}V_{11}& T_{11}V_{12}\\
J^{-1} (JT_{21})V_{11} &  J^{-1}(JT_{21})V_{12}\end{array}\right)\\
&=&  \left (\begin{array}{ccc}
T_{11}V_{11}  & T_{11}V_{12}\\
T_{21}V_{11}& T_{21}V_{12}\end{array}\right)\\
&=&TV.
\end{eqnarray*}.

Hene $\tilde{\phi}(T)= TV$, $T \in \mathcal{A}$ is an extension of $\phi$ and $\tilde{\phi} \in \mathcal{L(A)}.$
This completes the proof of (i) and (ii).

We come to the remaining part of the theorem. Let 
 %\begin {displaymath}
$U=  \phi(D).$  By Lemma \ref{third lemma},we have 
 \begin {displaymath}
U = \left (\begin{array}{ccc}
0 & 0\\
u_{21}& u_{22}\end{array}\right).\end{displaymath} .

%Now $D= D^2.$
%So $$U= \phi(D) = D\phi(D) =  \left (\begin{array}{ccc}
%0 & 0)\\
%u_{21}& u_{22}\end{array}\right)$$
Let $\phi_{U}: \mathcal{A}\rightarrow \mathcal{A}$ be given by 
$\phi_{U}(T)= TU.$ Then, $ \phi_{U}\in \mathcal {L(A)}.$
Let $\psi = \phi - \phi_{U_{|M_{2}}}.$
Then $ \psi \in  \mathcal{L}(M_{2},\mathcal{A})$.
Also, $\phi$ extends to $\tilde{\phi} \in  \mathcal {L(A)}$ if and only if $ \psi$ extends to some  $\tilde{\psi} \in \mathcal {L(A)}.$
We note that  \begin {displaymath}
\phi_{U} \left (\begin{array}{ccc}
T_{11} & 0\\
T_{21} & 0\end{array}\right) = 0 ,\end{displaymath} for $ T_{11} \in \mathcal{A}_{1}, T_{21}  \in \mathcal{A}_{5}.$ 
For $ t \in \mathcal{A}_{1}$, let \begin {displaymath}
L_{t}=  \left (\begin{array}{ccc}
t & 0 \\
0 & 0\end{array}\right) . \end{displaymath}

Then, for $ t,s \in \mathcal{A}_{1}$, $ L_{s}L_{t}= L_{st}.$
In particular, $A = L_{I_{X_{AH}}}.$
Let $ K \in \mathcal{A}_{3}.$
 Then $L_{K} \in M_{2}$, $\phi(L_{K}) = \psi(L_{K})$ 
and $L_{K} = AL_{K}.$
So $\phi(L_{K})= A\phi(L_{K}).$ 
This gives $ \phi(L_{K})$ has  the form \begin {displaymath}
 \left (\begin{array}{ccc}
S_{11} & S_{12}\\
0 & 0\end{array}\right) . \end{displaymath}.
Now, for $t \in \mathcal{A}_{1}, K \in \mathcal{A}_{3}$ we have
\begin {displaymath}
\phi(L_{tK})= \phi(L_{t}L_{K}) = L_{t}\phi(L_{K})= \left (\begin{array}{ccc}
tS_{11} & tS_{12}\\
0 & 0 \end{array}\right).\end{displaymath}.

Therefore the maps $\psi_{1}$ and $\psi_{2}$ on $\mathcal{A}_{3}$ to $\mathcal{A}_{1}$ and $\mathcal{A}_{2}$ respectively defined by $\psi_{1}(K)=S_{11}$,
$\psi_{2}(K)=S_{12}$ are bounded left $\mathcal{A}_{1}$-module homomorphisms.
Now, let $l \in \mathcal{A}_{5}.$ Then $Jl \in \mathcal{A}_{3}.
$
Also $ L^{l} =  \left (\begin{array}{ccc}
0 & 0\\
l & 0 \end{array}\right) \in M_{2}$ and $ \phi_{U}(L^{l}) =0.$
Now $AL^{l}= 0.$
So $A \phi(L^{l})=0.$
Therefore $\phi(L^{l})=\left (\begin{array}{ccc}
0 & 0\\
\beta_{21} & \beta_{22} \end{array}\right)$, for some $\beta_{21}, \beta_{22}.$ Now  $BL^{l}=L_{Jl}.$  
So 
\begin{eqnarray*}
B\left (\begin{array}{ccc}
0 & 0\\
 \beta_{21}& \beta_{22} \end{array}\right)\\ &= & \phi(L_{Jl})\\
&= &\left (\begin{array}{ccc}
\psi_{1}(Jl) & \psi_{2}(Jl)\\
0 & 0 \end{array}\right).\end{eqnarray*}

This gives, $J\beta_{21} = \psi_{1}(Jl)$, $J\beta_{22} = \psi_{2}(Jl).$
Consequently, Range $\psi_{1}(Jl)\subseteq Y$ and  Range $\psi_{2}(Jl)\subseteq Y$ , $\beta_{21}= J^{-1}\psi_{1}(Jl), \beta_{22}= J^{-1}\psi_{2}(Jl).$

Thus \begin{eqnarray*}
\phi(L^{l})&=& \psi(L^{l})\\
&=& \left (\begin{array}{ccc}
 0& 0\\
 J^{-1}\psi_{1}(Jl) & J^{-1}\psi_{2}(Jl) \end{array}\right) .\end{eqnarray*}
Hence $\phi(T)$ has the desired form and the rest follows from above.
 
\begin{theorem}\label{M2new}
$M_{2}$  cannot act as an interspersing submodule.
\end{theorem}
\pf We now utilise Remark 4.
Since $X_{AH}$ is not reflexive, there exists $\psi_{1}\in \mathcal{L}(\mathcal {A}_{3},\mathcal {A}_{1})$ such that $\psi_{1}$ is not extendible to $\mathcal{A}_{1}.$
Also by \cite{DS}, since $Y$ is not reflexive, there exists $\psi_{2} \in \mathcal{L}(\mathcal {A}_{3},\mathcal {A}_{2})$ such that $\psi_{2}$ is not extendible to $\mathcal{A}_{1}.$
Therefore by (ii) of Theorem \ref{m2} there exists $\phi \in  \mathcal{L}(M_{2},\mathcal{A})$ which is not extendible to $\mathcal{A}.$
\qed
\begin{corollary} 
Suppose $\it I$ is a closed left ideal of $\mathcal{A}= B(Z)$ containing $\mathcal{F}(Z).$ If $ \it I$ is contained in $M_{2}$ then $ \it I$ cannot serve as an interspersing Banach left $\mathcal{A}$-submodule. In particular, it is so for $M_1 \cap M_{2}.$ 
\end{corollary} 
\pf
By hypothesis,  $\it{I} \subseteq M_{2}.$ 
By Theorem \ref{M2new} there exists $\phi \in \mathcal{L}(M_{2},\mathcal{A})$ which is not extendible to $\mathcal{A}.$
Therefore $\phi_{| \it I} \in  \mathcal{L}(\it {I},\mathcal{A}) $ is also not extendible by Lemma 2.4.
\qed

\begin{rem} \cite{KL}\label{kl}
Let $\mathcal{T}_{2}$ be the algebra of upper triangular $2 \times 2$ matrices over the scalar field  $k.$
Then one can define $\phi:B(Z)\rightarrow \mathcal{T}_{2}$ by 
\begin{displaymath} 
T=\left (\begin{array}{ccc}
\alpha_{11}I_{X_{AH}}+K_{11}&\alpha_{12}J+K_{12}\\
K_{21} & \alpha_{22}I_{Y}+K_{22}\end{array}\right)\rightarrow \left (\begin{array}{ccc}
\alpha_{11}&\alpha_{12}\\
0 & \alpha_{22}\end{array}\right)
\end{displaymath}
and $\psi :\mathcal{T}_{2}\rightarrow B(Z)$ by 
\begin{displaymath} 
\left (\begin{array}{ccc}
\alpha_{11}&\alpha_{12}\\
0 & \alpha_{22}\end{array} \right) \rightarrow \left (\begin{array}{ccc}
\alpha_{11}I_{X_{AH}}&\alpha_{12}J\\
0 & \alpha_{22}I_{Y}\end{array}\right).\end{displaymath}
Clearly $ker(\phi)= K(Z)$ and $\phi_{0}\psi= I_{\mathcal{T}_{2}}.$
Also, it can be shown that for any left ideal $\mathcal{L}$ in $\mathcal{T}_{2}$, $ \mathcal{L} \rightarrow \phi^{-1}(\mathcal{L})$ is an order isomorphism 
of the lattice of left ideals of $\mathcal{T}_{2}$ onto the lattice of closed left ideals of $B(Z)$ that contain $K(Z).$ Apart from $ K(Z$), $M_{1}$, $M_{2}$ and $M_{1}\cap M_{2},$ the only such closed left ideals are of the form $ \mathcal{I}_{b}= \{a T_{b} + K : a \in k,~ K \in K(Z)\},$ where, for $b \in k$, $T_{b}= A + bB.$ Each of them is contained in $M_1.$ Arguments similar to those for Theorem \ref{m2} and \ref{M2new} give that there exist non-extendible bounded left $\mathcal A$-module homomorphisms on $\mathcal{I}_{b}$ to $B(Z).$ 
\end{rem}
\begin{theorem}
Every closed left ideal in $B(Z)$ other than $M_2,$ possesses an interspersing Banach left 
$\mathcal A$-submodule.
\end{theorem}
\pf We have only to combine Theorem \ref{dales}, Theorem \ref {banach algebra}, Theorem \ref{M1}, Theorem \ref{M2new}  and Remark \ref{kl}. 
\qed

\thanks { {\bf Acknowledgment}
 :We are grateful to D. Yost, T. Schlumprecht and  G. Godefroy for their valualble comments, suggestions and discussions. Ajit Iqbal Singh thanks the Indian National Science Academy for continuous support.}

\end{document}